\documentclass[a4paper,11pt]{amsart}

\usepackage{fullpage,amsxtra,amssymb,amsthm,amsmath,amscd}
\usepackage[latin1]{inputenc}

\renewcommand{\leq}{\leqslant}
\renewcommand{\geq}{\geqslant}

\newcommand{\sheaf}[1]{\mathcal{{#1}}}

\newcommand{\barre}[1]{\overline{{#1}}}

\newcommand{\Cc}{\mathbf{C}}

\newcommand{\Zz}{\mathbf{Z}}

\newcommand{\Rr}{\mathbf{R}}

\newcommand{\Qq}{\mathbf{Q}}
\newcommand{\Fp}{\mathbf{F}}

\newcommand{\mods}[1]{\,(\mathrm{mod}\,{#1})}

\newcommand{\ra}{\rightarrow}
\newcommand{\lra}{\longrightarrow}

\newcommand{\fleche}[1]{\stackrel{#1}{\lra}}

\DeclareMathOperator{\Spin}{Spin}

\DeclareMathOperator{\spec}{Spec}

\DeclareMathOperator{\rank}{rank}

\DeclareMathOperator{\ord}{ord}

\DeclareMathOperator{\frob}{Fr}
\DeclareMathOperator{\Gal}{Gal}

\DeclareMathOperator{\Tr}{Tr}

\DeclareMathOperator{\Aut}{Aut}

\DeclareMathSymbol{\gena}{\mathord}{letters}{"3C}
\DeclareMathSymbol{\genb}{\mathord}{letters}{"3E}

\newcommand{\consta}{\tfrac{1}{2N^2-2N+1}}

\theoremstyle{plain}
\newtheorem{theorem}{Theorem}

\newtheorem{lemma}[theorem]{Lemma}
\newtheorem{corollary}[theorem]{Corollary}

\newtheorem{proposition}[theorem]{Proposition}

\theoremstyle{remark}
\newtheorem{remark}[theorem]{Remark}

\theoremstyle{definition}

\begin{document}

\title{On the rank of quadratic twists of elliptic curvers over
  function fields}
\author{E. Kowalski}
\address{Universit\'e Bordeaux I - A2X\\
351, cours de la Lib\'eration\\
33405 Talence Cedex\\
France}
\email{emmanuel.kowalski@math.u-bordeaux1.fr}
\begin{abstract}
We prove quantitative upper bounds for the number of quadratic twists
of a given elliptic curve $E/\Fp_q(C)$ over a function field over a
finite field that have rank $\geq 2$, and for their average rank. The
main tools are constructions 
and results of Katz and uniform versions of the Chebotarev density
theorem for varieties over finite fields. Moreover, we conditionally
derive a bound in some cases where the degree of the conductor is
unbounded. 
\end{abstract}

\maketitle

Let first $E/\Qq$ be an elliptic curve over $\Qq$, and for fundamental
quadratic discriminants $d$, let $E_d$ denote the curve $E$ twisted by the
associated Kronecker character $\chi_d$. Goldfeld conjectured that
$E_d$ is most of the time of minimal rank compatible with the root
number of $E_d$, which in this case means
$$
\lim_{D\ra +\infty}{\frac{1}{|\{d\,\mid\, |d|\leq D\}|}
\sum_{|d|\leq D}{\rank E_d(\Qq)}}=\frac{1}{2}.
$$
\par
This conjecture has been refined
by Conrey, Keating, Rubinstein and Snaith~\cite{ckrs}, using
ideas based on Random Matrix Theory models and discretization
properties of special values of $L$-functions. For instance,
restricting the family to those $d$ for which the sign of the
functional equation of $L(E_d,s)$ is $+1$, they predict that for some
constants $c_E>0$ and $b_E\in\Rr$, we have
\begin{equation}\label{eq-ckrs}
|\{d\,\mid\, |d|\leq D,\text{ the root number of $E_d$ is $1$ and }\rank 
  E_d(\Qq)\geq 2\}| \sim c_E D^{3/4}(\log D)^{b_E}
\end{equation}
as $D\ra +\infty$, whereas the number of $d$ being counted is of size
$D$. 
\par
Note in particular that this predicts that there are few
curves with large rank, but gives also a lower bound for this
number. 
\par
From the analytic point of view, both conjectures are naturally
seen as statements about the order of vanishing of $L$-functions at
the central critical point, translated into the rank of $E_d$ by
assuming the Birch and Swinnerton-Dyer conjecture. This is indeed how
they arise, and one may expect to make progress on the analytic side
independently of the status of the Birch and Swinnerton-Dyer
conjecture. 
\par
The analogue of this problem for elliptic curves over function fields
has been developped by Katz~\cite{katz}: given an elliptic curve
$E/\Fp_q(t)$ (or over another function field over a finite field),
Katz shows how to construct various 
algebraic varieties $X/\Fp_q$ which parameterize twists of $E$ subject
to certain conditions. After a deep monodromy computation, he
obtains what can be considered as an analogue of Goldfeld's original
statement\footnote{\ Stated for the analytic rank, but from Ulmer's
  work on the Birch and Swinnerton-Dyer conjecture when the analytic
  rank is $\leq 1$, the case of algebraic rank follows, as explained
  in Section~\ref{sec-quad-twists}.}, up to the fact that the parameter
which gets large is not the conductor of the twisted curve, but rather
the degree of the coefficient field of the functions used to twist the
curve (see also~\cite[p. 134, 135]{ulmer} and
Section~\ref{sec-quad-twists}). 
\par
It is tempting to attack the more refined conjecture of~\cite{ckrs}
next, but there appear analytic difficulties in the application of Deligne's
equidistribution theorem. 
\par
We will show how to use a uniform version of Chebotarev's density
theorem (based on uniform estimates for $\ell$-adic Betti numbers
proved in~\cite{kow}) to
obtain a stronger quantitative form of the analogue of Goldfeld's
conjecture, under some monodromy assumptions which follow from the results
of Katz. This can be seen as a first progress towards the analogue of
the upper bound in~(\ref{eq-ckrs}) in this context; see
Corollary~\ref{cor-bounded-cond} and 
Proposition~\ref{pr-ortho} for precise statements (roughly, a small
power of $q$ is gained). 
\par
It would be very interesting also to obtain a lower
bound, but we do not consider this question; notice however that over
$\Qq$, fairly strong lower bounds for the occurence of algebraic rank
$\geq 2$ are known by
using results of sieve theory (see e.g.~\cite{gouvea-mazur}). This
method should presumably extend to the function field case (the lower
bound would probably be closer to the truth than our upper bounds
are), and lower bounds for occurence of algebraic rank $\geq 2$ are
also lower bounds for occurence of analytic rank $\geq 2$ in this
case. 
\par
In principle, we could  obtain results in situations where the
conductor is unbounded and the twists are restricted to special
one-parameter families. However, in that case we need rather stronger
monodromy results, and those do not seem available (on the other hand,
they are certainly within the realm of reason), even in special
cases. Still, we describe what could be true in the last section of
this paper.
\par
We should also mention that the results of Katz are in fact much more
general, and the method used here should adapt easily. We restrict our
attention to the case of quadratic twists of elliptic curves partly for
concreteness and partly in the hope of providing a reasonably readable
introduction to those remarkable results for readers with an analytic
number theory background. 
\par
\medskip
\textbf{Notation.} As usual, $|X|$ denotes the cardinality of a
set. 
By $f\ll g$ for $x\in X$, or $f=O(g)$ for $x\in X$, where $X$ is an
arbitrary set on which $f$ is defined, we mean synonymously that there
exists a constant $C\geq 0$ such that $|f(x)|\leq Cg(x)$ for all $x\in
X$. The ``implied constant'' is any admissible value of $C$. It may
depend on the set $X$ which is always specified or clear in context.
\par
\medskip
\textbf{Acknowledgments.} Work on this paper was prompted by
discussions by J. Keating. I wish also to thank N. Katz for explaining
a number of points concerning the behavior of monodromy for sheaves
with orthogonal symmetry. 

\section{A uniform version of the Chebotarev density theorem}

In this section we prove a general uniform Chebotarev density theorem
for varieties over finite fields. The main tools are the cohomological
methods and results developped notably by Grothendieck and Deligne; a
short and fairly concrete survey aimed at analytic number theorists,
which should be sufficient to explain the terminology and the proofs
below, can be found in~\cite[11.11]{ant}.
\par
We consider the following data: $U/\Fp_q$ is a smooth affine variety,
absolutely irreducible and of dimension $d\geq 1$ over a finite field of
characteristic $p$ with $q$ elements, $\ell\not=p$ is a prime number,
and $\rho\,:\, \pi_1(U,\barre{\eta})\ra G$ is a surjective 
map from the arithmetic fundamental group of $U$ (relative to the
geometric generic point of $U$) to a finite group $G$. 
We denote by
$$
G^g=\rho(\pi_1(\barre{U},\barre{\eta}))\subset G
$$
the image of the geometric fundamental
group of $U$, where $\barre{U}=U\times\barre{\Fp}_q$. Recall there are
exact sequences
\begin{equation}\label{eq-exact-seq}
\begin{CD}
1 @>>> \pi_1(\barre{U},\barre{\eta}) @>>> \pi_1(U,\barre{\eta})
@>d>> \hat{\Zz} @>>>  1\\
@. @VVV   @VVV  @V\varphi VV \\
1 @>>> G^g @>>>  G @>m>> \Gamma @>>>  1,
\end{CD}
\end{equation}
where the quotient $\Gamma$ thus defined is a finite cyclic
group. 
\par
For any $u\in U(\Fp_q)$, we denote by $\frob_u$ the geometric
Frobenius conjugacy class at $u$ in $\pi_1(U,\barre{\eta})$. In other
words, corresponding to 
the inclusion $\{u\}=\spec (\Fp_q)\ra U$, we have an induced homomorphism
$$
\pi_1(\Fp_q)=\Gal(\barre{\Fp}_q/\Fp_q)\ra \pi_1(U,\barre{\eta})
$$
and $\frob_u$ is the image, well-defined up to conjugation, of the
inverse of the generator $x\mapsto x^q$ of
$\Gal(\barre{\Fp}_q/\Fp_q)$. (If we consider $u$ as defined over a
bigger field, then $\frob_u$ changes, so $\frob_u$ is defined relative to
the field $\Fp_q$; it is often denoted $\frob_{u,\Fp_q}$ for this
reason, but we consider the base field as fixed in our statements). In
the exact sequence above, we have
$$
d(\frob_{u,q^n})=-n\in\hat{\Zz}.
$$
\par
The uniform Chebotarev density theorem is the following:

\begin{theorem}\label{th-chebo}
With notation as above, let $C\subset G$ be a
conjugacy-invariant subset such that $m(x)=\varphi(-1)$ for all $x\in
C$. Put
$$
\pi(C;q)=|\{u\in U(\Fp_q)\,\mid\, \rho(\frob_u)\in C\}|.
$$
\par
Assume that $G^g$ is of order prime to $p$. Then we have
\begin{equation}\label{eq-chebo}
\pi(C;q)=\frac{|C|}{|G^g|}|U(\Fp_q)|
+O(q^{d-1/2}|G|^{3/2}|C|^{1/2}),
\end{equation}
the implied constant depending only on
$\barre{U}=U\times\barre{\Fp}_q$. In particular, this holds uniformly 
with $q$ replaced by $q^n$ and $U$ by $U\times\Fp_{q^n}$,
$n\geq 1$.
\end{theorem}

\begin{proof}
This is essentially the same as the statement in
e.g.~\cite[Th. 4.1]{chavdarov}, except that we
have to take care of the uniformity. Let $f$ denote the
characteristic function of $C$ and let
$$
f(g)=\sum_{\pi}{\alpha(\pi)\Tr\pi(g)}
$$
be its Fourier expansion in terms of the irreducible representations
$\pi$ of $G$, which we realize as homomorphisms
$$
\pi\,:\, G\ra GL(\deg(\pi),E)
$$
for some finite extension $E/\Qq_{\ell}$, which can be chosen
independent of $\pi$. If $\pi=\psi$ is a character of $\Gamma$, i.e.,
$\pi$ is trivial on $G^g$, we have
$\Tr\pi(\rho(\frob_u))=\psi(\varphi(-1))$ by the assumption on $C$,
and
$$
\sum_{\psi}{\alpha(\psi)\Tr\pi(\rho(\frob_u))}
=\psi(\varphi(-1))\sum_{\psi}{\frac{1}{|G|}\sum_{x\in C}{\psi(x)}}
=\frac{|C|}{|G^g|}.
$$
\par
Now applying the Fourier expansion to $f(\rho(\frob_u))$ we find
therefore 
$$
\pi(C;q)=\frac{|C|}{|G^g|}|U(\Fp_q)|
+\sum_{\pi|G^g\not=1}{\alpha(\pi)
\sum_{u\in U(\Fp_q)}{\Tr\pi(\rho_{\ell}(\frob_u))}}.
$$
The inner sum is the sum of local traces for the 
representation 
$$
\pi\circ\rho\,:\, \pi_1(U,\barre{\eta})\ra GL(\deg(\pi),E),
$$
which can be seen as a lisse
$\barre{\Qq}_{\ell}$-adic sheaf, denoted $\pi(\rho)$. Since
the image of $\rho$, hence of $\pi(\rho)$, is finite,
this sheaf is pointwise pure of weight $0$.  
By the Grothendieck-Lefschetz trace formula we have
$$
\sum_{u\in U(\Fp_q)}{\Tr\pi(\rho(\frob_u))}=\sum_{i=0}^{2d}{
(-1)^i\Tr(\frob\,\mid\, H^i_c(\barre{U},\pi(\rho)))},
$$
where $\barre{U}=U\times\barre{\Fp}_q$. 
\par
In terms of the
geometric fundamental group $\pi_1(\barre{U},\barre{\eta})$, the
coinvariant description of $H^{2d}_c$ gives
$$
H^{2d}_c(\barre{U},\pi(\rho))
=E^{\deg(\pi)}_{\pi_1(\barre{U},\barre{\eta})}(-d)
=E_{G^g}^{\deg(\pi)}(-d)=0.
$$
since $\pi$, being non-trivial when restricted to $G^g$, can not
contain the trivial representation, simply because the space of
invariants under $G^g$ is a subrepresentation of $G$, which is
irreducible. 
\par
Moreover, by Deligne's Theorem, the eigenvalues of the geometric
Frobenius $\frob$ acting on each $H^i_c(\barre{U},\pi(\rho))$ are
algebraic integers with absolute value in $\Cc$ of modulus $\leq
q^{i/2}$. Thus we find
$$
\Bigl|\sum_{u\in U(\Fp_q)}{\Tr\pi(\rho(\frob_u))}\Bigr|
\leq q^{d-1/2}\sigma'_c(\barre{U},\pi(\rho)),
$$
where
$$
\sigma'_c(\barre{U},\pi(\rho))=\sum_{i<2d}{\dim
  H^i_c(\barre{U},\pi(\rho))}. 
$$
\par
It only remains to bound the quantity
$$
\sum_{\pi|G^g\not=1}{|\alpha(\pi)|\sigma'_c(\barre{U},\pi(\rho))}
$$
uniformly in terms of $\pi$ and $\rho$. By Proposition~3.6
of~\cite{kow}, using the assumption that $|G^g|$ is prime to
$p$, there exists a constant $C\geq 0$, depending only on $\barre{U}$,
such that 
$$
\sigma'_c(\barre{U},\pi(\rho))\leq C|G|(\deg\pi)
$$
for all $\pi$. Hence
$$
\sum_{\pi|G^g\not=1}{\sigma'_c(\barre{U},\pi(\rho))}
\leq C|G|\sum_{\pi}{|\alpha(\pi)|\deg\pi}
$$
and by Cauchy's inequality and the standard properties of
representations of finite groups we get
$$
\sum_{\pi}{\deg\pi}\leq \Bigl(\sum_{\pi}{|\alpha(\pi)|^2}\Bigr)^{1/2}
\Bigl(\sum_{\pi}{(\deg\pi)^2}\Bigr)^{1/2}
= \sqrt{|C|}\sqrt{|G|}.
$$
Putting these inequalities together yields the stated result.
\end{proof}

\begin{remark}
Since
$$
|U(\Fp_q)|=q^d+O(q^{d-1/2})
$$
by the Lang-Weil estimate, the implied constant depending only on
$\barre{U}$, we can also rephrase the result as
$$
\pi(C;q)=\frac{|C|}{|G^g|}q^d
+O(q^{d-1/2}|G|^{3/2}|C|^{1/2}).
$$
\end{remark}

Here is a variant of this theorem when the variety $U$ is a smooth
affine curve and the map $\rho$ arises by reduction from a
torsion-free $\Zz_{\ell}$-adic sheaf, using Proposition~3.1
of~\cite{kow} (or indeed, since we assume tameness, the last part of
Theorem~4.1 in~\cite{chavdarov}) instead of 
Proposition~3.6. In this case, the dependence on $U$ can be made
explicit (and the error term is improved) which allows certain
interesting applications (see the last section). 

\begin{theorem}\label{th-curve}
Let $U/\Fp_q$ be a smooth geometrically irreducible affine curve, realized
as an open dense subset of a smooth projective curve $C/\Fp_q$ of
genus $g$, with $m=|(C-U)(\barre{\Fp}_q)|$ ``points at infinity''. Let
$\ell\not=p$ be a prime number, let $\sheaf{F}$ be a tame torsion free
lisse $\Zz_{\ell}$-adic sheaf of rank 
$N$, and let $\barre{\sheaf{F}}=\sheaf{F}/\ell\sheaf{F}$ be its
reduction modulo $\ell$. Denote
$$
\rho_{\ell}\,:\, \pi_1(U,\barre{\eta})\ra GL(N,\Fp_{\ell})
$$
the corresponding continuous representation and put
$$
G_{\ell}=\rho_{\ell}(\pi_1(U,\barre{\eta})),\quad
G^g_{\ell}=\rho_{\ell}(\pi_1(\barre{U},\barre{\eta})),
$$
let $\Gamma_{\ell}$ be the quotient and $m$, $\varphi$ as
in~\emph{(\ref{eq-exact-seq})} in this case.
\par
Then for any conjugacy invariant subset $C(\ell)\subset G_{\ell}$ such 
that $m(C(\ell))=\varphi(-1)$,
$$
\pi(C(\ell);q)=\frac{|C(\ell)|}{|G_{\ell}|}|U(\Fp_q)|
+O((m+g)q^{1/2}|G_{\ell}|^{1/2}|C|^{1/2}),
$$
where the implied constant is \emph{absolute}.
\end{theorem}

\begin{proof}
As above we find by Fourier expansion that
$$
\pi(C(\ell);q)=\frac{|C(\ell)|}{|G_{\ell}|}|U(\Fp_q)|
+O(q^{1/2}\mathfrak{S})
$$
where 
$$
\mathfrak{S}=\sum_{\pi|G^g_{\ell}\not=1}{|\alpha(\pi)|
\sigma'_c(\barre{U},\pi(\rho_{\ell}))}
$$
and the implied constant is absolute (in fact, can be taken to be
equal to $1$).
\par
By Proposition~3.1 of~\cite{kow} (see eq.~(3.1) and note that the term
$w$ there vanishes by the tameness assumption), we have
$$
\sigma'_c(\barre{U},\pi(\rho_{\ell}))\leq (1-\chi_c(\barre{U},\Qq_{\ell}))
(\deg\pi)
$$
where
$$
\chi_c(\barre{U},\Qq_{\ell})=\dim H^0_c(\barre{U},\Qq_{\ell})-
\dim H^1_c(\barre{U},\Qq_{\ell})+\dim H^2_c(\barre{U},\Qq_{\ell})
$$
is the Euler-Poincaré characteristic of $\barre{U}$. This is equal to
$2-2g-m$ (because it is ``additive'', so equal to
$\chi_c(\barre{C})-\chi_c(\barre{C}-\barre{U})=2-2g-m$), hence
$$
\sigma'_c(\barre{U},\pi(\rho_{\ell}))\leq (2g+m-1)(\deg\pi).
$$
Summing over $\pi$ yields the estimate
$$
\mathfrak{S}\leq (m+2g-1)\sum_{\pi}{|\alpha(\pi)|\deg\pi}
\leq (m+2g-1)\sqrt{|G_{\ell}||C|},
$$
hence the result claimed.
\end{proof}

\begin{remark}
There are other uniform versions of the Chebotarev density theorem for
curves, for instance~\cite[Pr. 5.16]{fried-jarden}, which is written
and proved in a style closer to the standard number field
case. 
\end{remark}

\section{Application to sheaves with orthogonal monodromy}
\label{sec-ortho}

We will now apply the uniform Chebotarev density theorem to 
reductions of a lisse $\ell$-adic sheaf
with ``orthogonal'' symmetry. The method turns out to be essentially
identical with that used by Serre in his applications of the
Chebotarev density theorem over number fields
(see~\cite{serre-chebo}). 
\par
As in the previous section we start with a smooth affine absolutely
irreducible variety of dimension $d\geq 1$ defined over a finite field
$\Fp_q$ of characteristic $p$. Let $N\geq 1$ be an integer and
$\ell\not=p$ a prime number such that $GL(N,\Fp_{\ell})$ has order
prime to $p$. This also implies that for all $\nu\geq 1$, the group
$GL(N,\Zz/\ell^{\nu}\Zz)$ has order prime to $p$.
\par
Consider now a lisse integral torsion free $\ell$-adic sheaf
$\sheaf{F}_{\ell}$ of rank $N$
on $U$; equivalently, consider a continuous representation
$$
\tau_{\ell}\,:\, \pi_1(U,\barre{\eta})\ra GL(N,\Zz_{\ell}).
$$
We assume that 
$\sheaf{F}_{\ell}$ is equipped with a non-degenerate symmetric
pairing
$$
\langle \cdot,\cdot\rangle \,:\,
\sheaf{F}_{\ell}\otimes\sheaf{F}_{\ell}
\ra \Zz_{\ell},
$$
equivalently, $\tau_{\ell}$ acts on $\Zz_{\ell}^N$ by transformations
leaving invariant a non-degenerate symmetric pairing. We denote by
$O(N,\Zz_{\ell})$ the whole group of transformations leaving this
pairing invariant (which depends on the equivalence class of the
pairing, but this will not be of any importance).
\par
We make the following assumption of large monodromy on $\tau_{\ell}$:
\begin{equation}\label{eq-geom-arith}
[O(N,\Zz_{\ell}):\tau_{\ell}(\pi_1(U,\barre{\eta}))]<+\infty.
\end{equation}
\par
For any $\nu\geq 1$ we can consider the reduction modulo $\ell^{\nu}$
of $\sheaf{F}_{\ell}$, i.e. 
$\barre{\sheaf{F}}_{\nu}=\sheaf{F}_{\ell}/\ell^{\nu}\sheaf{F}_{\ell}$,
which corresponds to the maps
$$
\rho_{\nu}\,:\, \pi_1(U,\barre{\eta})\ra GL(N,\Zz/\ell^{\nu}\Zz),
$$
and we put
$$
G_{\nu}=\rho_{\nu}(\pi_1(U,\barre{\eta})).
$$
\par
The assumption~(\ref{eq-geom-arith}) ensures that the groups
$G_{\nu}$ are also ``large'' when $\nu$ is large enough : the index
\begin{equation}\label{eq-indices}
[O(N,\Zz/\ell^{\nu}\Zz):G_{\nu}]
\end{equation}
is \emph{bounded} for $\nu\geq 1$ (the reductions of any finite set of
coset representatives for $\tau_{\ell}(\pi_1(U,\barre{\eta}))$ in
$O(N,\Zz_{\ell})$ give coset representatives for $G_{\nu}$ in
$O(N,\Zz/\ell^{\nu}\Zz)$). 
\par
\medskip
We will apply the Chebotarev density theorem to sets $C_{\nu}\subset
G_{\nu}$ defined by reference with the forced eigenvalues that exist
for some orthogonal matrices. Precisely, recall that for $A\in O(N,k)$
(for an arbitrary field $k$ of odd characteristic and an arbitrary
non-degenerate symmetric bilinear form on $k^N$), the following 
``functional equation''  
\begin{equation}\label{eq-fn-eq-ortho}
T^NP(1/T)=\det(-A)P(T)
\end{equation}
holds for the polynomial $P=\det(1-AT)\in k[T]$. At $T=1$ this implies
that if $\det(-A)=(-1)^N\det A=-1$, we have $P(1)=0$, i.e., $1$ is
then an eigenvalue of $A$.
\par
If $N$ is even, that means that any matrix $A\in O(N,k)$ with
determinant $-1$ has eigenvalue $1$. If $N$ is odd, that means that
any orthogonal matrix in $SO(N,k)$ has eigenvalue $1$.  
We will say that $A$ has \emph{no extra vanishing} if $1$ is an
eigenvalue of $A$ with minimal multiplicity compatible with this. We
denote by $O^{ev}(N)$ the set of orthogonal matrices which \emph{have}
extra vanishing. We thus see that $O^{ev}(N)$ is an algebraic variety
defined over the same base as the orthogonal group under consideration,
given by $O^{ev}(N)=O_1\cup O_2$,
where
\begin{gather*}
O_1=\{A\in O(N)\,\mid\, \det(A)=-1\text{ and } \det(1-A)=0\}
\\
O_2=\{A\in SO(N)\,\mid\, \tfrac{d}{dT}\det(1-TA)\Bigl|_{T=1}=0\}
\end{gather*}
if $N$ is odd and
\begin{gather*}
O_1=\{A\in SO(N)\,\mid\, \det(1-A)=0\}
\\
O_2=\{A\in O(N)\,\mid\, \det(A)=-1\text{ and }
\tfrac{d}{dT}\det(1-TA)\Bigl|_{T=1}=0\}
\end{gather*}
if $N$ is even.
\par
We see that $O^{ev}(N)$ intersects each of the two connected
components of $O(N)$ in a closed hypersurface. 
\par
Denote by $C_{\nu}\subset G_{\nu}$ the image of
$$
O^{ev}(N,\Zz_{\ell})\cap \tau_{\ell}(\pi_1(U,\barre{\eta}))
$$
by reduction modulo $\ell^{\nu}$. We have the following simple lemma: 

\begin{lemma}\label{lm-cl}
\emph{(1)} We have
\begin{equation}\label{eq-gnu}
|G_{\nu}|\leq \ell^{\nu N(N-1)/2},
\end{equation}
for $\nu\geq 1$, the implied constant depending on $N$ and
$\ell$. 
\par
\emph{(2)} We have
\begin{equation}\label{eq-upper}
|C_{\nu}||G_{\nu}|^{-1}\ll \ell^{-\nu}
\end{equation}
for all $\nu\geq 1$, the implied constant depending on $N$, $\ell$ and
the index~\emph{(\ref{eq-geom-arith})}.
\end{lemma}

\begin{proof}
(1) The order of $G_{\nu}$ is at most the order of
$O(N,\Zz/\ell^{\nu}\Zz)$ so the upper bound follows easily. 
\par
(2) Because each $O^{ev}(N)$ is a hypersurface (i.e., of codimension
$1$ in each component of $O(N)$), the result follows from the bounds
$$
|C_{\nu}|\ll \ell^{\nu(\dim O(N)-1)},\quad
|G_{\nu}|\gg \ell^{\nu\dim O(N)},
$$
the first being a consequence of e.g.~\cite[Th. 8]{serre-chebo} and
the second of the boundedness of the indices~(\ref{eq-indices}) i.e.,
the finiteness of~(\ref{eq-geom-arith}).
(The upper bound can also be proved by more or less direct counting). 
\end{proof}

We now make the following
observation: for a point $u\in U(\Fp_q)$, the condition that 
$\frob_u$ has extra vanishing acting on $\sheaf{F}_{\ell}$ 
implies that $\rho_{\nu}(\frob_u)\in C_{\nu}$ for \emph{all} $\nu\geq
1$.  This leads to the basic bound
\begin{equation}\label{eq-basic}
|\{u\in U(\Fp_q)\,\mid\, \frob_u\text{ has extra vanishing on }
\sheaf{F}_{\ell}\}|
\leq \pi(C_{\nu};q)
\end{equation}
with notation as in Theorem~\ref{th-chebo}, which is valid for all
$\nu$. 
\par
We will use this to prove:

\begin{proposition}\label{pr-ortho}
With notation and assumptions as above, in particular under the
monodromy assumption~\emph{(\ref{eq-geom-arith})}, we have 
$$
|\{u\in U(\Fp_q)\,\mid\, \frob_u\text{ has extra vanishing on }
\sheaf{F}_{\ell}\}|
\ll q^{d-c}
$$
where $c=\consta$,
the implied constant depending on $\barre{U}$, $\ell$, $N$ and the
index~\emph{(\ref{eq-geom-arith})}. 
\end{proposition}

For instance, the gain is $c=1/13$ for $N=3$, $c=1/25$ for $N=4$.

\begin{proof}[Proof of Proposition~\ref{pr-ortho}]
Applying~(\ref{eq-basic}) and Theorem~\ref{th-chebo}, we derive
$$
|\{u\in U(\Fp_q)\,\mid\, \frob_u\text{ has extra vanishing}\}|
\leq \frac{|C_{\nu}|}{|G_{\nu}|}q^d+O(q^{d-1/2}|G_{\nu}|^{3/2}
|C_{\nu}|^{1/2})
$$
for any $\nu\geq 1$, with an implied constant depending only on
$\barre{U}$. 
\par
Using the bounds from the lemma, this leads to
$$
|\{u\in U(\Fp_q)\,\mid\, \frob_u\text{ has extra vanishing}\}|
\ll q^d\ell^{-\nu}+q^{d-1/2}\ell^{\nu (N(N-1)-1/2)},
$$
for $\nu\geq 1$ with an implied constant depending $\barre{U}$,
$\ell$, $N$ and the index~(\ref{eq-geom-arith}). For $q$ large enough
there exists $\nu$ such that
$$
\frac{1}{\ell}q^{c}\leq \ell^{\nu}\leq q^{c}\text{, with }
c=\consta,
$$
and then we have
$$
\ell^{-\nu}\leq \ell q^{-c},\text{ and }
q^{-1/2}\ell^{\nu (N(N-1)-1/2)}\leq q^{-c}
$$
so that taking this $\nu$ yields
$$
|\{u\in U(\Fp_q)\,\mid\, \frob_u\text{ has extra vanishing}\}|
\ll q^{d-c},
$$
the implied constant depending only on the parameters indicated
(the constant may need to be increased to account for small
values of $q$ where the $\ell$ above can not be found). 
\end{proof}


\begin{remark}
One could also approach the same problem using the large sieve method
of~\cite{kow}. However, because the ``sieving'' here is to detect an
algebraic condition, it turned out to be just as efficient using the
Chebotarev density theorem, and in fact using it for fixed $\ell$ as
above instead of requiring monodromy assumptions for many $\ell$. 
The large sieve would give a very slight improvement (roughly
corresponding to replacing $|C|$ in~(\ref{eq-chebo}) by the number of
conjugacy classes in $G$), but we avoid the necessary modification
necessary to deal with the fact that the orthogonal groups do not
satisfy the linear disjointness condition of loc. cit.
\end{remark}

We conclude by giving an equivalent rephrasing of the monodromy
assumption~(\ref{eq-geom-arith}).

\begin{proposition}\label{pr-large-mono}
Let $U/\Fp_q$ be a smooth absolutely irreducible variety over
$\Fp_q$. Let $\ell\not=p$ be a prime number and $\tau_{\ell}\,:\,
\pi_1(\barre{U},\barre{\eta})\ra 
GL(N,\Zz_{\ell})$ a continuous representation of the geometric
fundamental group of $U$. Assume $\rho_{\ell}$ takes value in
$O(N,\Zz_{\ell})$ for some non-degenerate bilinear form. If the
geometric monodromy group of $\tau_{\ell}$, i.e., the Zariski closure
of the image of 
$\tau_{\ell}$ in $GL(N,\barre{\Qq}_{\ell})$, contains
$SO(N)$, then we have
$$
[O(N,\Zz_{\ell}):\tau_{\ell}(\pi_1(U,\barre{\eta}))]<+\infty.
$$
\end{proposition}

In concrete terms, the group $G^g$ is the
set of matrices $x\in GL(N,\barre{\Qq}_{\ell})$ for which $f(x)=0$
whenever $f$ is a polynomial function on $GL(N,\barre{\Qq}_{\ell})$
(involving possibly $1/\det(x)$) that vanishes identically on
$\Gamma$, and the assumption is that $G^g\supset
SO(N,\barre{\Qq}_{\ell})$. 
\par
Thus, intuitively, both the condition~(\ref{eq-geom-arith})
and the assumption of the proposition are statements saying that
$\rho_{\ell}$ has a ``large'' image, and the proposition
shows that those two different meanings are in fact quite close.
\par
The condition~(\ref{eq-geom-arith}) may be called the ``old'' way of
stating this, whereas the assumption on the geometric monodromy group
is the more ``modern'' style; it is the usual language of the works of
Katz for instance; compare
with~\cite{serre-chebo}. In~\cite{serre-p-div}, Serre
attributes the shift to Grothendieck. 

\begin{proof}
Since the image of $\pi_1(U,\barre{\eta})$ is larger than that of
$\pi_1(\barre{U},\barre{\eta})$, it suffices to show that
$$
[O(N,\Zz_{\ell}):\tau_{\ell}(\pi_1(\barre{U},\barre{\eta}))]<+\infty.
$$
\par
Let $\Gamma$ denote the image of the geometric fundamental group. 
In the $\ell$-adic topology induced from $GL(N,\Qq_{\ell})$,
$\Gamma$ is compact (by continuity of $\tau_{\ell}$, as the
fundamental group is compact) inside the compact group
$O(N,\Zz_{\ell})$. The point is that from
$G^g\supset SO(N)$, it follows that $\Gamma\cap SO(N,\Zz_{\ell})$ is
also open in $SO(N,\Zz_{\ell})$, still for the $\ell$-adic
topology (see~\cite[Cor. p. 120]{serre-p-div} for instance). But then
the finiteness of $[SO(N,\Zz_{\ell}):\Gamma\cap SO(N,\Zz_{\ell})]$
is immediate from the existence of Haar measure $\mu$ with total mass one on
$SO(N,\Zz_{\ell})$, since $\Gamma$ as an open set must satisfy
$\mu(\Gamma)>0$. In fact,
$$
[SO(N,\Zz_{\ell}):\Gamma\cap SO(N,\Zz_{\ell})]=\frac{1}{\mu(\Gamma)}.
$$
\par
The desired finiteness of $[O(N,\Zz_{\ell}):\Gamma]$ is obviously a trivial
consequence of this.
\end{proof}

\section{Twists of elliptic curves over function fields}
\label{sec-quad-twists}

We now explain how the result of the previous section apply
to the study of extra vanishing for families of twists of
elliptic curves over function fields.
\par
We first survey the construction by Katz of varieties parameterizing
twists of elliptic curves over function fields over finite fields
(see~\cite[Intro., V]{katz} and also~\cite[2,6,7]{ulmer}). 
\par
We assume for simplicity that the characteristic $p$ is not $2$ or
$3$. 
Let $C/\Fp_q$ be a smooth projective curve of genus $g$, absolutely
irreducible, $K=\Fp_q(C)$ the function field of $C$ and $E/K$ an
elliptic curve, which is assumed to have non-constant
$j$-invariant. If $C=\mathbf{P}^1$ is the projective line, the twists
considered can be described concretely as follows: take a Weierstrass
equation for $E$ of the type
$$
y^2=x^3+a(t)x^2+b(t)x+c(t)
$$
with $a$, $b$, $c\in\Fp_q(t)$, and let $f\in \Fp_q[t]$ be a
(squarefree) polynomial. Then the twist $E_f$ of $E$ by $f$ is the
elliptic curver $E_f/K$ with equation
$$
f(t)y^2=x^3+a(t)x^2+b(t)x+c(t),
$$
and coefficients of $f$ serve as ``algebraic'' parameters for twists. 
\par
In the greater generality described, $f$ is chosen to be a rational
function on $C$ with a prescribed set of poles $D$ (an effective
divisor on $C$) and with $\deg D$ distinct zeroes, none of which is a
place of bad reduction of $E$. Katz shows that if $\deg(D)\geq 2g+1$,
this set of functions is the set of $\Fp_q$-rational points of a
smooth geometrically connected algebraic variety $X/\Fp_q$  (which
depends on $D$ and $E$). 
\par
Now for any $f\in X(\Fp_{q^n})$, $n\geq 1$,
there is a twist $E_f$ defined over $\Fp_{q^n}(C)$, generalizing the
above description. For any prime $\ell\not=p$, the $L$-functions of
all the twisted curves can be encoded in the ``local'' behavior
of a certain lisse $\ell$-adic sheaf $\sheaf{T}_{\ell}$ on $X$, of
rank $N$ which is independent of $\ell$, 
corresponding to a representation 
$$
\tilde{\tau}_{\ell}\,:\, \pi_1(U,\barre{\eta})\ra GL(N,\Zz_{\ell}).
$$
\par
Precisely, for any
rational point $f\in X(\Fp_q)$, we have the identity
\begin{equation}\label{eq-l-function}
L(E_f/K,T)=\det(1-T\frob_x\,\mid\, \sheaf{T}_{\ell}),\text{ or }
L(E_f/K,s)=\det(1-q^{-s}\frob_x\,\mid\, \sheaf{T}_{\ell}).
\end{equation}
\par
This sheaf is constructed by Katz~\cite[Ch. V]{katz}. Moreover, Katz
shows, as consequences of general properties of 
étale cohomology, that $\sheaf{T}_{\ell}$ is punctually pure of weight $2$,
and there exists a natural non-degenerate symmetric pairing
$$
\sheaf{T}_{\ell}\otimes\sheaf{T}_{\ell}\ra \Zz_{\ell}(-2).
$$
\par
This means that the image of $\pi_1(U,\barre{\eta})$ by $\tilde{\tau}_{\ell}$
is contained in the group $CO(N,\Zz_{\ell})$ of similitudes for this
pairing, and the image of $\frob_{f,q^n}$, for $u\in X(\Fp_{q^n})$, is
a similitude with ``multiplier''
$m(\tilde{\tau}_{\ell}(\frob_{u,q^n}))=q^{2n}$. In addition, the image of the
geometric fundamental group is contained in the group
$O(N,\Zz_{\ell})$ of orthogonal transformations for the pairing. Of
course, we have an exact sequence
$$
1\ra O(N,\Zz_{\ell})\ra CO(N,\Zz_{\ell})\fleche{m} \Zz_{\ell}^{\times}
\ra 1.
$$
\par
The $L$-functions~(\ref{eq-l-function}) have central critical point at
$s=1$, i.e., at $T=q^{-1}$, and it is convenient to make a Tate twist
to translate it to $s=0$, i.e., $T=1$. 
\par
So we consider the sheaf  $\sheaf{T}_{\ell}(1)$ 
instead of $\sheaf{T}_{\ell}$, which corresponds to taking the
representation
$$
\tau_{\ell}(x)=\tilde{\tau}_{\ell}(x)q^{-d(x)}
$$
of $\pi_1(U,\barre{\eta})$, where $d$ is the ``degree'' map
in~(\ref{eq-exact-seq}). The twisted representation
$\tau_{\ell}$ coincides with $\tilde{\tau}_{\ell}$ on
$\pi_1(\barre{U},\barre{\eta})$ (by~(\ref{eq-exact-seq})). On the
other hand, since $\pi_1(U,\barre{\eta})$ is topologically generated by
the $\frob_{f,q^n}$, $f\in X(\Fp_{q^n})$, $n\geq 1$, and since
$$
m(q^{-d(\frob_{u,q^n})}\tilde{\tau}_{\ell}(\frob_{u,q^n}))=
q^{-2n}q^{2n}=1,
$$
it follows that $\tau_{\ell}(\pi_1(U,\barre{\eta}))\subset
O(N,\Zz_{\ell})$.
\par
With $\sheaf{F}_{\ell}=\sheaf{T}_{\ell}(1)$, this provides us with all
the data occuring in 
Section~\ref{sec-ortho}. Now the point is that Katz has shown by a
deep monodromy computation that the
condition~(\ref{eq-geom-arith}) holds for $D$ suitably chosen. We
state a precise version: 

\begin{proposition}[Katz]
If the divisor $D$  satisfies the conditions
\begin{equation}\label{eq-cond-div}
\deg(D)\geq 4g+4,\quad 2g-2+\deg(D)\geq \max(144, 2s),
\end{equation}
where $s$ is the number of places of bad reduction of $E/K$, then we
have
$$
[O(N,\Zz_{\ell}):\tau_{\ell}(\pi_1(U,\barre{\eta}))]<+\infty.
$$
\end{proposition}

\begin{proof}
Under the condition stated, Katz has shown (see~\cite[p. 15]{katz}
for a summary) that the geometric monodromy group $G^g$ associated to
$\sheaf{T}_{\ell}(1)$ contains $SO(N)$. Therefore, we can apply
Proposition~\ref{pr-large-mono}.
\end{proof}

\begin{corollary}\label{cor-bounded-cond}
Let $C/\Fp_q$ be a smooth absolutely irreducible projective curve of
genus $g\geq 0$, 
$E/\Fp_q(C)$ an elliptic curve with non-constant $j$-invariant, $D$ an
effective divisor on $E$ of degree $\geq 2g+1$, and $X=X(D,E)$ the
associated parameter space 
for twists. Assume that $p>N+2$ and that the twisting sheaves
$\sheaf{T}_{\ell}(1)$ satisfy~\emph{(\ref{eq-geom-arith})};
for instance assume that~\emph{(\ref{eq-cond-div})} holds.
\par
Then, for any $n\geq 1$, the number $V_n$ of
twisting parameters $f\in X(\Fp_{q^n})$ such 
that the $L$-function of $E_f/\Fp_{q^n}(C)$ vanishes at $s=1$ with order
strictly larger than that imposed by the functional equation satisfies
\begin{equation}\label{eq-nn}
V_n\ll q^{n(\dim X-c)}
\end{equation}
with $c=\consta$, the implied constant depending on $D$,
$E$, and $p$. In particular this set has density $0$ as
$n\ra +\infty$.
\end{corollary}

\begin{proof} 
The condition $p>N+2$ implies that there exists a non-zero congruence
class $a$ modulo $p$ such that $\ell\equiv a\mods{p}$ implies that
$GL(N,\Fp_{\ell})$ is of order prime to $p$ (see e.g.~\cite[Lemma
7.5.1]{katz}). Pick such a prime $\ell$, and then apply 
Proposition~\ref{pr-ortho} and~(\ref{eq-l-function}) to the sheaf
$\sheaf{T}_{\ell}(1)$.
\end{proof}

It is interesting to notice that, together with the work of Ulmer
(see~\cite{ulmer-2}) on the Birch and Swinnerton-Dyer conjecture which
shows that  
$$
\rank E_f(\Fp_{q^n}(C))=\ord_{s=1}{L(E_f/\Fp_{q^n}(C),s)}
$$
\emph{if} the right-hand side is at most $1$, this also gives a strong
version of the analogue of Goldfeld's Conjecture for the \emph{algebraic}
rank. To state it, we assume for simplicity that the image of the geometric
fundamental group by the representation $\tau_{\ell}$ is not contained
in $SO(N,\Zz_{\ell})$. Katz has shown (\cite[Ex. 8.3.4.1]{katz})
that this is the case for instance if $E/K$ has multiplicative
reduction at a point $s\in S$.

\begin{proposition}\label{pr-goldfeld}
With assumptions as in the previous corollary, assume moreover
that the image of the geometric
fundamental group by the representation $\tau_{\ell}$ is not contained
in $SO(N,\Zz_{\ell})$. Then
we have
$$
\sum_{f\in X(\Fp_{q^n})}{\rank E_f(\Fp_{q^n}(C))}
=\frac{1}{2}|X(\Fp_{q^n})|+O(q^{n(\dim X-c)})
$$
for $n\geq 1$, the implied constant depending on $D$, $E$ and $p$.
\end{proposition}

\begin{proof}
The point is that all the twists have ``analytic rank''
bounded by the rank $N$ 
of $\sheaf{T}_{\ell}$, since their $L$-functions are all polynomials of degree
$N$. Since for an elliptic curve $E/\Fp_{q^n}(C)$  there is the a-priori
inequality (due to Tate~\cite{tate})
$$
\rank E(\Fp_{q^n}(C))\leq \ord_{s=1} L(E/\Fp_{q^n}(C),s),
$$
the twists also have algebraic rank bounded by $N$,
so the contribution to the average rank of those few $f$ for which the
analytic rank of $E_f$ is $\geq 2$ is small. Precisely, notice that
$$
\rank E_f(\Fp_{q^n}(C))=\frac{1-W(E_f)}{2}+\beta(E_f)
$$
where $W(E_f)$ is the sign of the functional equation and $\beta(E_f)$
has the property that $\beta(E_f)=0$ if the 
analytic rank of $E_f$ is $\leq 1$ (this is the result of Ulmer), and
$|\beta(E_f)|$ is bounded ($\leq N+1$) for all $f$ and $n$. Thus the
sum of $\beta(E_f)$ is 
$$
\ll q^{n(\dim X-c)}
$$
for $n\geq 1$ by~(\ref{eq-nn}). On the other hand
$$
\sum_{f\in X(\Fp_{q^n})}{\frac{1-W(E_f)}{2}}
=\frac{|X(\Fp_{q^n})|}{2}
-\frac{1}{2} \sum_{f\in X(\Fp_{q^n})}{W(E_f)}.
$$
and we have
$$
W(E_f)=(-1)^N\det(\tau_{\ell}(\frob_f))
$$
where $\tau_{\ell}$ is as before the representation which
corresponds to $\sheaf{T}_{\ell}$ (this is
simply~(\ref{eq-fn-eq-ortho}) and~(\ref{eq-l-function})). 
\par
The assumption that the image of $\tau_{\ell}$ is not inside $SO$
implies that the character $\det\circ \tau_{\ell}$ of order $2$ is
non-trivial on the geometric fundamental group. 
Thus $H^{2d}_c(\barre{X},\det\circ\tau_{\ell})=0$ and by the Riemann
hypothesis we have 
$$
\Bigl|\sum_{f\in X(\Fp_{q^n})}{W(E_f)}\Bigr|\leq
q^{n(\dim X-1/2)}\sum_{0\leq i<2\dim X}{\dim
  H^i_c(\barre{X},\det\circ\tau_{\ell})}
\ll q^{n(\dim X-1/2)}
$$
for $n\geq 1$, finishing the proof.
\end{proof}

\section{Examples of one parameter families of twists}

This section is a concrete illustration of the previous
section. We will use some of the intermediate statements proved by Katz
in~\cite{katz} to restrict our attention to one-parameter families of
twists (i.e., replace the big parameter space $X$ of the previous
section by a curve) where the analog of
Corollary~\ref{cor-bounded-cond} and Proposition~\ref{pr-goldfeld}
still hold. (In fact, it is by finding such one-parameter families
inside the larger spaces that Katz shows that the
geometric monodromy groups for those parameter spaces contain $SO(N)$).
\par
To simplify, we will only consider curves over 
$\mathbf{P}^1$ and twists by polynomials; this allows us to write down
explicit equations (in other words, $C=\mathbf{P}^1$, $g=0$ and the
divisor $D$ is $d(\infty)$ in the notation of the previous section).
\par
We start with an elliptic curve given by a fairly general Weierstrass
equation (assuming always that the characteristic $p$ is $\geq 5$)
$$
E\,:\,y^2=x^3+a(t)x^2+b(t)x+c(t)
$$
with polynomials $a$, $b$ and $c\in \Fp_q[t]$. Denote by $S$ the set
of points in $\mathbf{A}^1$ where $E$ has bad reduction. 
\par
In Chapter~5 of~\cite{katz}, two different types of one-parameter
families with ``large'' monodromy are described. We will consider
special cases of the first one (Theorem~5.4.1 of loc. cit.).
Adapting the
terminology found there and 
in~\cite[Ch. 2]{katz}, we say that a polynomial $f\in \barre{\Fp}_q[t]$
of degree $d\geq 1$ is \emph{of Lefschetz type} if the following
conditions hold: 
\par
\indent \quad (i) $f$ has $d$ distinct zeros in $\barre{\Fp}_q$; 
\par
\indent \quad (ii) $f'$ has $d-1$ distinct zeros, and those have distinct
images by $ f$. 
\par
Moreover we say that $f$ is of \emph{Katz-Lefschetz type} for $E$ if
$f$ is of Lefschetz type, and 
\par
\indent \quad (i') no two points of $S$ have the same image by $f$;
\par
\indent \quad (ii') $f(s)\not=0$ for $s\in S$; 
\par
\indent \quad (iii') the fibers $f^{-1}(f(s))$ all  have  $d$ distinct
elements for $s\in S$.
\par
Fix a polynomial $f$ of Katz-Lefschetz type. Let $V_f$ denote the
(finite) variety, defined over $\Fp_{q}$, of critical values for $f$,
i.e.
$$
V_f=f(\{x\,\mid\, f'(x)=0\})\cup f(S),
$$
and let $U_f=\mathbf{A}_1-V_f$. 
Then Theorem~5.4.1
of~\cite{katz} states that if $d\geq \max(146,2|S|)$,
the one-parameter family of quadratic
twists with equations  
$$
E_{f,\alpha}\,:\, (f(t)-\alpha)y^2=x^3+a(t)x^2+b(t)x+c(t),
$$
with $\alpha\in U_f$ 
has associated $\ell$-adic twisting sheaves $\sheaf{F}_{\ell,f}$ on
$U_f$, of fixed rank $N_f$, such that the geometric monodromy
group (on $U_f$) of $\sheaf{F}_{\ell,f}$ contains $SO(N_f)$  for
$\ell\not=p$. 
In particular this sheaf satisfies the condition~(\ref{eq-geom-arith}).

\begin{corollary}
Let $E/\Fp_q(t)$ be an elliptic curve with non-constant $j$-invariant
and at least one finite place of multiplicative reduction.
Let $f$ be a polynomial of Katz-Lefschetz type for $E$ of degree
$d\geq \max(146,2|S|)$.
\par
\emph{(i)} For $n\geq 1$, the
number $V_n$ of $\alpha\in U_f(\Fp_{q^n})$ where
$E_{f,\alpha}/\Fp_{q^n}(t)$ has extra vanishing 
satisfies
$$
V_n\ll q^{n(1-c)}
$$
with $c=\tfrac{1}{2N_f^2-2N_f+1}$, the implied constant depending on
$E$, $f$ and $p$.
\par
\emph{(ii)} If the image of the geometric fundamental group is not
contained in $SO(N,\Zz_{\ell})$, we have
$$
\sum_{\alpha\in U_f(\Fp_{q^n})}{\rank E_{f,\alpha}(\Fp_{q^n}(t))}
=\frac{1}{2}|U_f(\Fp_{q^n})|+O(q^{n(1-c)}),
$$
for $n\geq 1$, the implied constant depending on
$E$, $f$ and $p$.
\end{corollary}

We could in fact state a slightly better result using the Chebotarev density
theorem for curves instead of the general version (see the next
section).
\par
To be completely concrete, we will now take a specific example. Let
$E/\Fp_p(t)$ be the following variant of the Legendre elliptic curve:
\begin{equation}\label{eq-legendre}
E\,:\, y^2=x(x+1)(x-t).
\end{equation}
\par
Note that $E$ has \emph{multiplicative} reduction at $\infty$ and at
the points in $S=\{0,-1\}$. 
\par
Now consider the following polynomials:
$$
f_d=t^d-dt-1\in \Fp_{p}[t].
$$

\begin{lemma}
\emph{(i)} If $p\nmid d(d-1)$ and $(p-1,d-1)=1$, then $f_d$ is of
Lefschetz type. 
\par
\emph{(ii)} If in addition $p\nmid d+1$, then $f_d$
is of Katz-Lefschetz type for the above curve $E/\Fp_p(t)$.
\end{lemma}

\begin{proof}
The derivative of $f_d$ is $f'_d=d(t^{n-1}-1)$ so since $p\nmid d$,
the roots of $f'_d$ are the $(d-1)$-st roots of unity. Since $p\nmid
d-1$, there are $d-1$ of them in $\barre{\Fp}_p$. Now for $\mu$ a zero
of $f'_d$ we have
$$
f_d(\mu)=\mu(1-d)-1.
$$
This already shows that the values of $f_d$ at zeros of $f'_d$ are
distinct. 
\par
Notice that this formula also shows $f_d(\mu)\not=0$ because otherwise
$\mu$ would be in the prime field $\Fp_p$, so that $\mu=1$ by the
assumption $(p-1,d-1)=1$, the equation becomes $1-d=d$, but again
$p\nmid d$ excludes this case. 
Since $f_d^{-1}(x)$, for $x\in \barre{\Fp}_p$, has
$d$ elements except if $x$ is in the set $\{f_d(\mu)\}$, it follows
that $f_d$ has $d$ distinct roots in $\barre{\Fp}_p$. Altogether, this
establishes the first assertion that $f_d$ is of Lefschetz type.
\par
For the second, the conditions $(p-1,d-1)=1$ and $p\geq 5$ imply that
$d$ is even. We compute $f_d$ at the points in $S$: we have
$f_d(0)=-1$ and $f_d(-1)=(-1)^d+d-1=d$. So $f_d(0)\not=
f_d(-1)$ since $p\nmid d+1$. Moreover $f_d(0)=-1$ is not of the form
$f_d(\mu)=\mu(1-d)-1$ as above, since $p\nmid d$. Similarly
$f_d(-1)=d$ is not of this form: $d=\mu(1-d)-1$ implies again that
$\mu$ is in the prime 
field, so $\mu=1$, and again $p\nmid d$ shows that $d+1=1-d$ is
impossible. So neither $0$ nor $-1$ is in a fiber over a zero of
$f'_d$, which means that the fibers over points of $S$ contain $d$
distinct elements.
\end{proof}

\begin{remark}
So for $p=5$, we have found explicit polynomials of Katz-Lefschetz type
for $E$ of any even degree $d$ with $d\equiv 2\mods{5}$, i.e.,
$d\equiv 2\mods{10}$. 
\par
In general, the density of integers $d\geq 1$ satisfying the
conditions of the lemma is
$$
\frac{\varphi(p-1)}{p}\Bigl(1-\frac{3}{p}\Bigr)>0
$$
for all $p\geq 5$.
\end{remark}

Let $d$ be any integer satisfying the condition of the lemma. We then
have the one-parameter family of twists
$$
E_{d,\alpha}\,:\, (f_d(t)-\alpha)y^2=x(x+1)(x-t),
$$
or equivalently (change $y$ to $(f_d-\alpha)y$)
$$
E_{d,\alpha}\,:\, y^2=(t^d-dt-1-\alpha)x(x+1)(x-t)
$$
over $\Fp_p(t)$, with parameter $\alpha$ in the complement $U_d$ of the
finite variety of critical values for $f_d$, which has
$d-1+2=d+1$ points defined over $\barre{\Fp}_{p}$. 
\par
Let $\sheaf{F}_{d,\ell}$ denote the twisting sheaf
$\sheaf{T}_{\ell}(1)$ for this subfamily.
As observed by Katz~\cite[Lemma 7.5.1]{katz}, the twist sheaves
associated to quadratic twists of elliptic curves are always tame in
characteristic $p\geq 5$, so $\sheaf{F}_{d,\ell}$ is tame.
The rank $N_d$ of
$\sheaf{F}_{d,\ell}$ is computed in~\cite[Lemma 5.1.3, p. 16]{katz}
and is given  
by
$$
N_d=2d
$$
for $d$ even. In particular, this means that the (degree of the)
conductor of the twists goes to infinity when $d\ra +\infty$.
\par
As a special case of Theorem~5.4.1 of~\cite{katz}, if $d\geq 146$, the
geometric monodromy group for $\sheaf{F}_{d,\ell}$ is the full
orthogonal group $O(N_d)$ for all $\ell$.
\par
So specializing again the previous corollary we get:
\begin{corollary}
Let $p\geq 5$ be prime, let $d\geq 146$, $f_d$ and $U_d$ be as
above. For $n\geq 1$, the number $V_n$ of $\alpha\in U_d(\Fp_{q^n})$ for
which the twisted Legendre curve
$$
y^2=(t^d-dt+1-\alpha)x(x+1)(x-t)
$$
over $\Fp_{p^n}(t)$ has extra vanishing satisfies
\begin{equation}\label{eq-bound-legendre}
V_n\ll p^{n(1-c)},
\end{equation}
with $c=\tfrac{1}{2(d^2+1)}$, the implied constant depending on $d$
and $p$. 
\end{corollary}

\section{Twists with unbounded conductor}
\label{sec-unbounded}

This section is speculative. The idea
is to exploit the strong bound~(\ref{eq-bound-legendre})
to prove a variant of Proposition~\ref{pr-goldfeld} for a
family of twists more closely resembling the quadratic twists of
elliptic curves over $\Qq$, namely one where the conductor (i.e.,
essentially, in this case, the degree of the $L$-function) increases,
so that the rank of the elliptic curves is not uniformly
bounded.  The speculation consists in
the fact that the result obtained is conditional on monodromy
assumptions which are stronger than currently known.
\par
We still work with the curve~(\ref{eq-legendre}) of the previous
section for concreteness. Take a sequence of polynomials $f_n$ of
Katz-Lefschetz type with increasing degrees $d_n$. We will use a
simple subscript $n$ for all the objects of the last section which
would otherwise require to be subscripted by either $f_n$ or
$d_n$. For instance, we denote by $\sheaf{F}_{n,\ell}$ the twisting
sheaf for the $1$-parameter family corresponding to $f_n$ and let
$N_n$ denote its rank. 
\par
We make the
following strong assumption 
\begin{gather}
\text{For all $n$ and \emph{all} odd $\ell\not=p$,
the image of the representation}\nonumber\\
\rho_{n,\ell}\,:\, \pi_1(U_{d_n},\bar{\eta})\ra 
GL(N_n,\Fp_{\ell})\label{eq-strong-mono}
\\
\text{corresponding to the reduction 
$\sheaf{F}_{n,\ell}/\ell\sheaf{F}_{n,\ell}$  is of \emph{bounded
index} in $O(N_n,\Fp_{\ell})$}.\nonumber
\end{gather}

See the final paragraphs of the paper for
comments on the plausibility of this.\footnote{.\ Note also
  that~(\ref{eq-strong-mono}) could be
replaced without much change by a ``vertical'' version, namely that
for some \emph{fixed} $\ell\not=p$, the index of the image of
$\pi_1(U_{d_n},\barre{\eta})\ra O(N_n,\Zz/\ell^{\nu}\Zz)$ is bounded
for all $n\geq 1$ and $\nu\geq 1$.}
\par
We denote by $G_{n,\ell}$ the image of $\rho_{n,\ell}$
and by $B$ a bound for its index in $O(N,\Fp_{\ell})$ valid for all
$n$ and $\ell\not=p$. 
\par
Corresponding to Lemma~\ref{lm-cl} we need the following uniform
version for $\nu=1$, which we make a little bit more precise:

\begin{lemma}\label{lm-cl2}
\emph{(1)} For all odd primes $\ell$ and all $n\geq 1$ we have
$$
|G_{n,\ell}|\leq \ell^{N_n(N_n-1)/2},
$$
the implied constant depending only on the bound $B$ for the index of
$G_{n,\ell}$ in $O(N_n,\Fp_{\ell})$.
\par
\emph{(2)} Let $C_{n,\ell}$ be the set of $g\in G_{n,\ell}$ with extra
vanishing. We have
$$
\frac{|C_{n,\ell}|}{|G_{n,\ell}|}\ll \frac{1}{\ell}
$$
for all odd primes $\ell$ and $n\geq 1$, the implied constant
depending only on $B$, provided that $\ell\geq N_n^2$.
\end{lemma}

\begin{proof}
(1) The size of $G_{d,\ell}$ is bounded by that of $O(N_n,\Fp_{\ell})$
for which the existing formulas immediately give the result
stated. 
\par
(2) We bound $C_{n,\ell}$ by the number of elements with extra
vanishing in $O(N_n,\Fp_{\ell})$. In general, for $O(N,\Fp_{\ell})$,
the latter (say $R(N)$) is written as follows: 
$$
R(N)=\sum_{g}{|\{ x\in O(N,\Fp_{\ell})\,\mid\, \det(1-Tx)=g\}|}
$$
where $g$ runs over characteristic polynomials of elements of
$O(N,\Fp_{\ell})$ which have extra vanishing. It is clear that the
number of possible $g$ is $\leq \ell^{N-1}$. For each $g$, we count
the inner quantity by the same method as
in~\cite[Proof of Th. 3.5]{chavdarov} (with adaptations necessary
because the orthogonal group is not simply connected like the
symplectic group) which shows 
that it is  
$$
\ll (\ell+1)^{N(N-1)/2}(\ell-1)^{-N}
$$
(with absolute implied constant),
so we get
$$
R(N)\leq \ell^{-1} \Bigl(\frac{\ell}{\ell-1}\Bigr)^N
(\ell+1)^{N(N-1)/2}
$$
and because $|O(N,\Fp_{\ell})|\geq (\ell-1)^{N(N-1)/2}$, this yields
the result after an application of the mean value theorem.
\end{proof}

Here is the hypothetical result with unbounded conductors.

\begin{proposition}
Let $p\geq 5$ be prime, let $d_n$ for $n\geq 1$ be an
increasing sequence of integers such that each $d=d_n$ satisfies
$p\nmid d(d-1)(d+1)$ and $(p-1,d-1)=1$, and such that
$d_n^3\leq n$ for $n\geq 1$. Assume the monodromy
hypothesis~\emph{(\ref{eq-strong-mono})} for
the sequence of Katz-Lefschetz polynomials $f_n=t^{d_n}-d_nt+1$.
\par
\emph{(i)} We have for $n\geq 1$
$$
|\{\alpha\in U_{d_n}(\Fp_{p^n})\,\mid\, 
E_{d_n,\alpha}\text{ has extra vanishing}\}|\ll
n^{1/3}p^{n-\tfrac{1}{4}n^{1/3}}.
$$
\par 
\emph{(ii)} We have for $n\geq 1$
$$
\sum_{\alpha\in U_{d_n}(\Fp_{p^n})}{
\rank E_{d_n,\alpha}(\Fp_{p^n}(t))
}
=\frac{p^n}{2}+O(n^{2/3}p^{n-\tfrac{1}{4}n^{1/3}}).
$$
\end{proposition}

\begin{proof}
Because of all the assumptions, the
Chebotarev density theorem for curves (Theorem~\ref{th-curve}, with
$g=0$, $m=d_n+1$) and Lemma~\ref{lm-cl2} imply that for $n\geq 1$ and
$\ell\not=p$ we have
$$
|\{\alpha\in U_{d_n}(\Fp_{q^n})\,\mid\, 
E_{d_n,\alpha}\text{ has extra vanishing}\}|\ll
\frac{p^n}{\ell}+d_np^{n/2}\ell^{A}
$$
with 
$$
A=\frac{N_n(N_n-1)}{2}-\frac{1}{2},
$$
the implied constant depending only on $B$ if $\ell\geq
N_n^2=4d_n^2$. 
\par
Since $2(A+1)=N_n(N_n-1)+1\leq 4d_n^2$, taking $\ell$ between
$p^{n/2(A+1)}$ and $2p^{n/2(A+1)}$ gives
$$
|\{\alpha\in U_{d_n}(\Fp_{p^n})\,\mid\, 
E_{d_n,\alpha}\text{ has extra vanishing}\}|\ll
d_np^{n-\tfrac{n}{4d_n^2}}
$$
for $n\geq 1$, with an implied constant depending only on $B$; the
assumption $d_n^3\leq n$ ensures that for $n$ large enough ($n\geq 6$
suffices) we have 
$$
n\log p\geq 2d_n^2(\log 4d_n^2)
$$
hence $\ell\geq 4d_n^2$, and the implied constant can be raised to
absord the values $n\leq 5$, as well as those for which $d_n<146$.
\par
For part (ii), the reasoning 
is as in the proof of Proposition~\ref{pr-goldfeld}, using Ulmer's
result about the Birch and Swinnerton-Dyer conjecture. The
contribution of the twists with analytic rank $\geq 2$ is estimated
using~(i) and the trivial bound
$$
\rank  E_{d_n,\alpha}(\Fp_{p^n}(t))\leq N_n = 2d_n\leq 2n^{1/3}
$$
(so there too the uniformity of our estimates in terms of $n$ is -- or
would be! -- important). 
\end{proof}

\begin{remark}
In terms of the parameter $X=p^n$, the error terms have the following
shape:
$$
n^{1/3}p^{n-\tfrac{1}{4}n^{1/3}}\ll (\log X)^{1/3}X
\exp(-\tfrac{1}{4}(\log X)^{1/3})
$$
which may look more familiar to analytic number theorists.
\end{remark}

We finish by commenting on our monodromy
assumption~(\ref{eq-strong-mono}). First of
all, for \emph{fixed} $n$, the uniformity in terms of $\ell$ is part
of the standard conjectures (see e.g.~\cite[10.3?,10.7?]{serre-conj})
about the variation of images of $\ell$-adic representations. 
\par
In addition, since Katz has 
shown that the ``rational'' geometric monodromy group is always equal to
$O(N)$, it is a consequence of a result of
Larsen~\cite[Th. 3.17]{larsen} that 
for a set of primes $\ell$ of density $1$, the geometric monodromy group
modulo $\ell$ contains the image in $O(N,\Fp_{\ell})$ of the spin
group $\Spin(N,\Fp_{\ell})$, which is of index $2$ in
$SO(N,\Fp_{\ell})$ and $4$ in $O(N,\Fp_{\ell})$ (this complication
arises because 
$O(N,\Fp_{\ell})$ is neither connected nor simply connected).
Larsen's result is quite difficult
(it uses the classification of simple finite groups), and the set of
primes it produces is not easy to control.
\par
Another example, still for fixed rank, is a result proved by Gabber
concerning the monodromy of Kloosterman sheaves which is explained
in~\cite[Ch. 12]{katz-gsksm}. Roughly speaking, the integral monodromy
group associated to families of Kloosterman sums in an even number $n$
of variables is ``big'' for all $\ell$ large enough, depending on $n$,
but again not in an easy way to describe for varying $n$.
\par
When the rank is increasing, it is in fact not clear if the uniform
bound we postulate is 
coherent with the general philosophy concerning $\ell$-adic
representations. The reason is that this variation of $n$ does not
fall into a well-understood theoretical framework: the ``family'' we
consider is one only inasmuch as we manage to deal with its individual
terms and get similar results; this is much the same as the case of
``families'' of classical automorphic $L$-functions, for which
convincing examples exist abundantly without an \emph{a priori}
definition. Still, in the case of elliptic curves $E/\Qq$, a similar
result is expected: recall that Serre showed that for a fixed $E$
without CM, there
exists $L_E$ such that for any prime $\ell>L_E$, the map
$$
\rho_{E,\ell}\,:\, \Gal(\barre{\Qq}/\Qq)\ra 
\Aut(E[\ell])\simeq GL(2,\Fp_{\ell}) 
$$
is surjective (i.e., the fields obtained by adjoining to $\Qq$ the
coordinates of the $\ell$-torsion points of $E$ are ``as big as
possible''). Then the conjectured statement that such an $L$ as above
exists which ``works'' for all elliptic curves $E/\Qq$ (without CM) can
be seen as an analogue of our assumption (see e.g.~\cite[Question~2,
p. 199]{serre-chebo}). It has been confirmed by Duke~\cite{duke-nex}
that this can be done for ``almost all'' curves.
\par
Finally, we can turn for encouragement to at least one similar situation
where a result of the desired type has been unconditionally proved.  
Let $f\in \Fp_q[x]$ be a fixed polynomial of
degree $2g$ with $2g$ distinct roots in $\barre{\Fp}_q$, and consider
the family of hyperelliptic curves of genus $g\geq 1$ with equations
$$
C_{\alpha}\,:\,y^2=f(x)(x-\alpha),
$$
over the open set $\mathbf{A}^1-f^{-1}(0)$, with projection
$\pi(x,y,\alpha)=\alpha$. The sheaves $R^1\pi_!\Fp_{\ell}$, which are
of rank $2g$ and 
admit symplectic symmetry, are used to
``globalize'' the family of $L$-functions of $C_{\alpha}$ (modulo
$\ell$). Jiu-Kang Yu has shown that the geometric monodromy group
is equal to $Sp(2g,\Fp_{\ell})$ for all $f$ and all $\ell\not=2$. 
This is one of the main examples in~\cite{chavdarov}. It is also used
to give some results uniform in $g$ in~\cite{kow} (which are in fact 
of a rather more delicate nature).

\end{document}